\documentclass[12pt,reqno]{amsart}
\usepackage{amssymb,verbatim,enumerate,ifthen}
\usepackage[mathscr]{eucal}
\usepackage[utf8]{inputenc}
\usepackage[T1]{fontenc}
\def\N{\mathbb{N}}
\def\R{\mathbb{R}}
\def\Q{\mathbb{Q}}
\def\Z{\mathbb{Z}}
\def\F{\mathbb{F}}
\def\T{\mathbb{T}}

\def\conv{\mathop{\mbox{\rm conv}}\nolimits}

\newtheorem{theorem}{Theorem}
\newtheorem*{theorem*}{Theorem}
\def\Thm#1#2{\ifthenelse{\equal{#1}{*}}{\begin{theorem*}#2\end{theorem*}}
             {\begin{theorem}\label{T#1}#2\end{theorem}}}
\newtheorem{Atheorem}{Theorem}

\def\thm#1{Theorem~\ref{T#1}}

\newtheorem{proposition}[theorem]{Proposition}
\newtheorem*{proposition*}{Proposition}
\def\Prp#1#2{\ifthenelse{\equal{#1}{*}}{\begin{proposition*}#2\end{proposition*}}
             {\begin{proposition}\label{P#1}#2\end{proposition}}}
\def\prp#1{Proposition~\ref{P#1}}

\newtheorem{corollary}[theorem]{Corollary}
\newtheorem*{corollary*}{Corollary}
\def\Cor#1#2{\ifthenelse{\equal{#1}{*}}{\begin{corollary*}#2\end{corollary*}}
             {\begin{corollary}\label{C#1}#2\end{corollary}}}
\def\cor#1{Corollary~\ref{C#1}}

\newtheorem{lemma}[theorem]{Lemma}
\newtheorem*{lemma*}{Lemma}
\def\Lem#1#2{\ifthenelse{\equal{#1}{*}}{\begin{lemma*}#2\end{lemma*}}
             {\begin{lemma}\label{L#1}#2\end{lemma}}}
\def\lem#1{Lemma~\ref{L#1}}

\newtheorem{remark}[theorem]{Remark}
\newtheorem*{remark*}{Remark}
\def\Rem#1#2{\ifthenelse{\equal{#1}{*}}{\begin{remark*}\rm #2\end{remark*}}
             {\begin{remark}\label{R#1}\rm #2\end{remark}}}
\def\rem#1{Remark~\ref{R#1}}

\newtheorem{example}[theorem]{Example}
\newtheorem*{example*}{Example}
\def\Exa#1#2{\ifthenelse{\equal{#1}{*}}{\begin{example*}\rm #2\end{example*}}
             {\begin{example}\label{Ex#1}\rm #2\end{example}}}

\def\eq#1{{\rm(\ref{E#1})}}
\def\Eq#1#2{\ifthenelse{\equal{#1}{*}}
  {\begin{equation*}\begin{aligned}[]#2\end{aligned}\end{equation*}}
  {\begin{equation}\begin{aligned}[]\label{E#1}#2\end{aligned}\end{equation}}}

\begin{document}
\begin{flushright}
\textit{Submitted to: Semigroup Forum} 
\end{flushright}
\vspace{5mm}

\date{\today}

\title[Convexity in abelian semigroup setting]
{Convexity and a Stone-type theorem for convex sets in abelian semigroup setting}

\author[W. Jarczyk]{Witold Jarczyk}
\address{Faculty of Mathematics, Informatics and Econometrics, 
University of Zielona G\'ora, \newline PL-65-516 Zielona G\'ora, Poland}
\email{w.jarczyk@wmie.uz.zgora.pl}

\author[Zs. P\'ales]{Zsolt P\'ales}
\address{Institute of Mathematics, University of Debrecen, 
H-4010 Debrecen, Pf.\ 12, Hungary}
\email{pales@science.unideb.hu}

\subjclass[2000]{Primary}

\thanks{
This research of the second author was realized in the frames of T\'AMOP 4.2.4. A/2-11-1-2012-0001 
”National Excellence Program – Elaborating and operating an inland student and researcher personal
support system”. The project was subsidized by the European Union and co-financed by the European 
Social Fund. This research of the second author was also supported by the Hungarian Scientific 
Research Fund (OTKA) Grant NK 81402.}

\begin{abstract}
In this paper, two parallel notions of convexity of sets are introduced in the abelian semigroup setting.
The connection of these notions to algebraic and to set-theoretic operations is investigated. A formula for 
the computation of the convex hull is derived. Finally, a Stone-type separation theorem for disjoint convex 
sets is established.
\end{abstract}

\maketitle

\section{Introduction}

Given a linear space $X$ and two disjoint convex sets $A_0,B_0 \subseteq X$, Stone's celebrated 
theorem asserts that there exist two disjoint convex sets $A,B \subseteq X$ such that
\Eq{*}{
A_0 \subseteq A, \qquad B_0 \subseteq B, \qquad \mbox{and} \qquad A \cup B = X
}
(see \cite{Sto46}). In other words, two disjoint convex sets can always be separated by two disjoint 
complementary convex sets. This basic result can be used to derive the standard Hahn--Banach-type theorems 
on the separation of convex sets and convex functions. For the details of this approach we refer 
to the book of Holmes \cite{Hol72}.

A completely analogous result was established by P\'ales \cite{Pal89d} for the separation of disjoint subsemigroups of an abelian semigroup $S$: Given two disjoint subsemigroups $A_0,B_0 \subseteq S$, 
there exist two disjoint subsemigroups $A,B \subseteq S$ such that
\Eq{*}{
A_0 \subseteq A, \qquad B_0 \subseteq B, \qquad \mbox{and} \qquad A \cup B = S.
}
This result is used to derive, for instance, the characterization theorem of quasideviation means
(cf.\ \cite{Pal82a}, \cite{Pal89b}).

The notion of convexity is generalized in various ways in the context of functions to the setting of abelian groups (cf.\ \cite{Jar10b}, \cite{JarLac09}, \cite{JarLac10}). In this paper we consider two concepts of convexity of subsets of any abelian semigroup and, as one of our main results, we derive a Stone-type separation theorem in that context.

Let $(S,+)$ denote an abelian semigroup throughout this paper. Given a subset $A \subseteq S$ 
and  an $n \in \N$, the sets $nA$, $n^{-1}A$, $[n]A$ are defined by
\Eq{nA}{
nA:=\{nx\mid x\in A\},\qquad n^{-1}A:=\{x\mid nx\in A\}, \\
[n]A := \{x_1+ \cdots +x_n \mid x_1, \ldots , x_n \in A\},\qquad
}
respectively. Obviously, the inclusions
\Eq{I1}{
  n(n^{-1}A)\subseteq A\subseteq n^{-1}(nA),\qquad A\subseteq n^{-1}([n]A), \qquad nA\subseteq[n]A
}
hold for every $n\in\N$ and $A\subseteq S$.

For a fixed $n \in \N$, a subset $A\subseteq S$ is called {\it $n$-convex} and {\it $n$-konvex}, if
\Eq{I2}{
  n^{-1}([n]A) \subseteq A  \qquad\mbox{and}\qquad [n]A \subseteq nA
}
hold, respectively. In fact, in view of the last two inclusions in \eq{I1}, a set $A$ is 
$n$-convex and $n$-konvex if and only if, in the respective case, there is equality in \eq{I2}.

The authors of the paper \cite{JarLac10}, while studying convex functions in a group setting, deal also 
with 2-convex sets and 2-konvex sets. 2-convex sets are called {\it convex} there, whereas 2-konvex subsets $A$ of the group are simply indicated by the equality $A+A=2A$. 

Observe that if $S$ is the additive group of a vector space over the field $\Q$, then the notions of $n$-convexity and $n$-konvexity coincide and a subset of $S$ is $n$-convex for all $n\in\N$ if and only if it is closed under rational convex combination, i.e., it is $\Q$-convex in the standard sense. 

More generally, if the semigroup $S$ is divisible by $n$, then $n$-convexity implies $n$-konvexity. In the case when this $n$-divisibility of $S$ is unique, both the notions of convexity coincide. In general, however, 
they are different whenever $n\geq 2$. The additive group $\mathbb Z$ of integers serves as an example of 
an $n$-convex set which is not $n$-konvex. On the other hand, consider the circle group $\T=\R/\Z$ 
identified with $[0,1)$, where addition is meant modulo 1, and let $A=[0,1/n]$. Then $[n]A=nA=[0,1)$ 
and $n\cdot(3/4)\in [0,1) =[n]A$, although $3/4 \not \in A$, and thus $A$ is $n$-konvex but not $n$-convex.

Let $\F$ be a nonvoid subset of $\N$. A set $A \subseteq S$ is said to be $\F$-{\it convex} 
(resp. $\F$-{\it konvex})  if it is $n$-convex (resp. $n$-konvex) for all $n \in \F$. If $A$ is $\N$-convex 
(resp. ${\mathbb N}$-konvex), then it is called {\it convex} (resp. {\it konvex}). Observe that the semigroup 
$S$ is automatically convex, however, it may not be konvex. On the other hand, for every $x\in S$, the 
singleton $\{x\}$ is konvex and, in general, it is not convex.

In Section 2 we compare both the notions and examine their algebraic and set-theoretical properties. Suitable 
examples show that some of them fail. It turns out that the notions of $n$-convexity and $n$-konvexity are, in a sense, complementary: a number of the properties is adhered only to one of the notions. In particular, the family of $n$-convex sets is closed under the intersection, but the family of $n$-konvex sets is not in general. For that reason, for every subset $A$ of the semigroup and for any nonempty set $\F$ of positive integers, we may consider only $\F$-convex hull $\conv_{\F}(A)$ of $A$. Given a multiplicative subsemigroup $\F$ of $\N$ we find the form of the set $\conv_{\F}(A)$. We describe also the $\N$-convex hull of the union of finitely many sets, in particular obtaining a kind of the drop theorem. The equivalence relation, 
determined by the partition of the semigroup into $\F$-convex hulls of singletons, is also studied.

Section 3 provides the notion of $\F$-disjointness and a a Stone-type theorem for the separation of 
$\N$-disjoint sets by complementary convex subsets of the semigroup. We give two proofs of it. The first one follows the previous results of the paper, whereas the second one makes use of the Stone-type theorem for separation of subsemigroups proved in \cite{Pal89d}.

The present paper deals with various aspects of the notions of convex subsets and konvex subsets of an abelian semigroup. The next step is to discuss possible notions of convex functions defined on a subset of a semigroup. Such a research in a group setting was presented in the paper \cite{JarLac10}; some further questions were answered in \cite{Pal82a} and \cite{JarLac09}.

\section{Properties of convex and konvex sets}

In addition to the inclusions in \eq{I1}, the following lemma summarizes the basic properties 
of the operations introduced in \eq{nA}.

\Lem{1}{For all $k, n \in \N$, and for all subsets $A,B\subseteq S$,
\Eq{L1}{
(k+n)A&\subseteq kA+nA, & (kn)A&=k(nA), & n(A+B)&=nA+nB,\\
[k+n]A&=[k]A+[n]A, & [kn]A&=[k]([n]A),& [n](A+B)&=[n]A+[n]B,\\
kA+nA&\supseteq (kn)(k^{-1}\!A+n^{-1}\!A), & (kn)^{-1}\!A&=k^{-1}(n^{-1}\!A), 
    & n^{-1}(A+B)&\supseteq n^{-1}\!A+n^{-1}\!B,\\ 
[k](n^{-1}\!A)&\subseteq n^{-1}([k]A), & [k](nA)&=n([k]A), & k(n^{-1}A)&\subseteq n^{-1}(kA).
}
The inclusions, in general, cannot be replaced by equalities. 
Provided that $A\subseteq B$, we also have
\Eq{L2}{
  nA\subseteq nB, \qquad n^{-1}A\subseteq n^{-1}B, \qquad [n]A\subseteq [n]B.
}
If $\{A_\gamma\mid\gamma\in\Gamma\}$ is a family of subsets of $S$, then
\Eq{L3}{
  n\bigg(\bigcap_{\gamma\in\Gamma} A_\gamma\bigg) &\subseteq \bigcap_{\gamma\in\Gamma} nA_\gamma,
  &\qquad
  n\bigg(\bigcup_{\gamma\in\Gamma} A_\gamma\bigg) &\supseteq \bigcup_{\gamma\in\Gamma} nA_\gamma,\\
  n^{-1}\bigg(\bigcap_{\gamma\in\Gamma} A_\gamma\bigg) &\subseteq \bigcap_{\gamma\in\Gamma} n^{-1}A_\gamma,
  &\qquad
  n^{-1}\bigg(\bigcup_{\gamma\in\Gamma} A_\gamma\bigg) &\supseteq \bigcup_{\gamma\in\Gamma} n^{-1}A_\gamma,\\
  [n]\bigg(\bigcap_{\gamma\in\Gamma} A_\gamma\bigg) &\subseteq \bigcap_{\gamma\in\Gamma} [n]A_\gamma,
  &\qquad
  [n]\bigg(\bigcup_{\gamma\in\Gamma} A_\gamma\bigg) &\supseteq \bigcup_{\gamma\in\Gamma} [n]A_\gamma.  
}}

\begin{proof} 
The proof of \eq{L1} follows from the associativity and commutativity of the semigroup operation in an 
elementary manner. \eq{L2} is straightforward and \eq{L3} directly follows from \eq{L2}.

In order to see that the five inclusions in \eq{L1} are proper in general, take the additive group $\Z$ and 
let $A=B$ be the set of odd numbers and let $k=n=2$. 
\end{proof}

For any subset $A\subseteq S$, define the sets
\Eq{CK}{
  C_A:=\{n\in\N\mid A \mbox{ is $n$-convex}\}\qquad\mbox{and}\qquad
  K_A:=\{n\in\N\mid A \mbox{ is $n$-konvex}\}.
}
The structural properties of these sets are contained in the following proposition.

\Prp{0}{For all subsets $A\subseteq S$, the set $K_A$ is a multiplicative subsemigroup of $\N$. 
On the other hand, if $n\in C_A$ and $k$ is a divisor of $n$, then $k\in C_A$.}

\begin{proof} Let $A$ be simultaneously $k$- and $n$-konvex. Then, using \lem{1}, we obtain
\Eq{*}{
  [kn]A=[k]([n]A)\subseteq [k](nA)=n([k]A)\subseteq n(kA)=(kn)A,
} 
proving that $A$ is $kn$-konvex. Therefore, $K_A$ is a multiplicative subsemigroup of $\N$.

Now let $n\in C_A$ and $k|n$. Let $m:=n/k$. To prove the $k$-convexity of $A$, choose $x\in k^{-1}([k]A)$. Then 
$kx\in [k]A$ and hence $nx=m(kx)\in m([k]A)\subseteq [mk]A=[n](A)$, i.e., $x\in n^{-1}([n]A)\subseteq A$.
\end{proof}

The next result deals with algebraic properties of the classes of $\F$-convex and $\F$-konvex sets.

\Prp{1}{Let $\F$ be a nonempty subset of $\N$. The family of $\F$-konvex sets is closed 
under the algebraic addition of sets and, for all $k\in\N$, under the multiplication by $k$ and under the operation $[k]$. The family of $\F$-konvex sets is not closed under the operation $k^{-1}$ in general.
For all $k\in\N$, the family of $\F$-convex sets is closed under the operation $k^{-1}$.
The family of $\F$-convex sets is not closed under the algebraic addition, under the multiplication by $k$ 
and under the operation $[k]$ in general.}

\begin{proof} Let $A$ and $B$ be $\F$-konvex sets and let $n\in\F$. Using \lem{1} and the inclusions 
$[n]A\subseteq nA$ and $[n]B\subseteq nB$, we get
\Eq{*}{
  [n](A+B)=[n]A+[n]B\subseteq nA+nB=n(A+B),
}
showing that $A+B$ is $n$-konvex for all $n\in\F$, hence it is $\F$-konvex, too.

Assuming that $k\in\N$ and that $A$ is $\F$-konvex, for all $n\in\F$, we have $[n]A\subseteq nA$, 
whence \lem{1} yields
\Eq{*}{
  [n](kA)&=k([n]A)\subseteq k(nA)=(kn)A=(nk)A=n(kA), \\
  [n]([k]A)&=[nk]A=[kn]A=[k]([n]A)\subseteq [k](nA)=n([k]A),
}
which proves that $kA$ and also $[k]A$ is $n$-konvex for all $n\in\F$.

Let $A$ be an $\F$-convex set and let $k\in\N$. Then, for all $n\in\F$, by \lem{1} and by using 
$n^{-1}([n]A)\subseteq A$, we 
have
\Eq{*}{
 n^{-1}\big([n](k^{-1}A)\big)\subseteq n^{-1}\big(k^{-1}([n]A)\big)=(nk)^{-1}([n]A)
  = (kn)^{-1}([n]A) = k^{-1}\big(n^{-1}([n]A)\big) \subseteq k^{-1}A,
}
proving that $k^{-1}A$ is $n$-convex for all $n\in\F$.

In what follows, let $\F$ be the singleton $\{n\}$ with $n>1$.
To show that the family of $n$-konvex sets is not closed under the operation $k^{-1}$ in general,
let $S:=\T=\R/\Z$. Then $\{0\}$ is $n$-konvex (because it is a singleton). On the other hand, for $k>1$, 
the set $k^{-1}\{0\}=\big\{0,\frac1k,\dots,\frac{k-1}{k}\big\}$ is not $n$-konvex whenever $n$ and $k$ are 
not relative primes. Indeed, then
\Eq{*}{
[n](k^{-1}\{0\})=[n]\Big\{0,\frac1k,\dots,\frac{k-1}{k}\Big\}
    =\Big\{0,\frac1k,\dots,\frac{k-1}{k}\Big\}=k^{-1}\{0\},
}
which has exactly $k$ elements. On the other hand, the set $n(k^{-1}\{0\})$ contains 
$\frac{k}{(n,k)}<k$ elements, therefore $[n](k^{-1}\{0\})$ cannot be equal to $n(k^{-1}\{0\})$,
proving that $k^{-1}\{0\}$ is not $n$-konvex.

To prove that the family of $n$-convex sets is not closed under the algebraic addition in general, let $S:= \Z$ and take $A:=\{0, n-1\}$ and $B:=\{0, 2n-1\}$. Then
\Eq{*}{
  [n]A=\{0, n-1, 2(n-1), \dots, n(n-1)\}
}
and
\Eq{*}{
  [n]B=\{0, 2n-1, 2(2n-1), \dots, n(2n-1)\}.
}
Since $n$ and $n-1$ are relative primes, it follows that $A$ is $n$-convex. Similarly, $B$ is $n$-convex as $n$ and $2n-1$ are relatively prime. Notice also that $A+B=\{0, n-1, 2n-1, 3n-2\}$. If $n=2k$ with some $k \in \N$, then we have
\Eq{*}{
  \frac1n\left(k(n-1)+k(2n-1)\right)=\frac1n\left((k+k)(n-1)+kn\right)=n-1+k,
}
and if $n=2k+1$  with some $k \in \N$, then we get
\Eq{*}{
  \frac1n\left((k+1)(n-1)+k(2n-1)\right)=\frac1n\left(((k+1)+k)(n-1)+kn\right)=n-1+k.
}
In both cases we get an element of the set $n^{-1}\left([n](A+B)\right)$ which does not belong to $A+B$. Therefore the set $A+B$ is not $n$-convex.

The set $\{0,1\}$ is $n$-convex in the additive group of $\Z$ since
\Eq{*}{
  n^{-1}([n]\{0,1\})=n^{-1}\{0,\dots,n\}=\{0,1\}.
}
On the other hand, $k\{0,1\}$ is $\{0,k\}$ and
\Eq{*}{
  n^{-1}([n]\{0,k\})=n^{-1}\{0, k, \dots, nk\},
}
which is strictly bigger than $\{0,k\}$ provided that $n$ and $k$ are not relative 
primes. Hence, in this case, $\{0,k\}$ is not $n$-convex.

Finally we show that the family of $n$-convex sets is not closed under the operation $[k]$ in general. To the aim let $S$ be the semigroup of nonnegative integers  different from $1$ and $A:=\{0,2\}$. Then $[n]A=\{0, 2, 4, \dots, 2n\}$ and $n^{-1}\big([n]A\big)=\{0,2\}$, and thus $A$ is $n$-convex. On the other hand, by \lem{1}, we have
\Eq{*}{
  n^{-1}\big([n]\big([k]A\big)\big)& =n^{-1}([nk]A)=n^{-1}\{0, 2, 4, \dots, 2nk\}\\
	& =\{0, 2, 3, 4, \dots, 2k\}\not\supseteq \{0, 2, 4, \dots, 2k\}=[k]A,
}
which means that $[k]A$ is not $n$-convex.
\end{proof}

\Rem{1}{In the case when the semigroup $S$ is uniquely divisible by $n$ then the notions of $n$-convexity and 
$n$-konvexity coincide, therefore, the family of $n$-convex sets is closed under the algebraic addition of 
sets and, for all $k\in\N$, under the multiplication by $k$ and under the operations $k^{-1}$, $[k]$.}

\Prp{2}{Let $\F$ be a nonempty subset of $\N$. Then the intersection of any family and the union 
of any chain of $\F$-convex sets is again $\F$-convex. In addition, the union of any chain of $\F$-konvex 
sets is again $\F$-konvex. The family of $\F$-konvex sets is not closed under the intersection in general.
Furthermore, the intersection of an $\F$-convex and an $\F$-konvex set is always $\F$-konvex. In particular, 
if $S$ is $\F$-konvex, then every $\F$-convex subset of $S$ is also $\F$-konvex.}

\begin{proof} Assume that $\{A_\gamma\mid\gamma\in\Gamma\}$ is a family of $\F$-convex sets.
Then, for any $n\in\F$, by \lem{1} and the $n$-convexity of the sets $A_\gamma$, we have
\Eq{*}{
  n^{-1}\bigg([n]\bigcap_{\gamma\in\Gamma}A_\gamma\bigg)
  \subseteq n^{-1}\bigg(\bigcap_{\gamma\in\Gamma}[n]A_\gamma\bigg)
  \subseteq \bigcap_{\gamma\in\Gamma}n^{-1}([n]A_\gamma)
  \subseteq \bigcap_{\gamma\in\Gamma}A_\gamma,
}
which demonstrates that the intersection of any family of $n$-convex sets is again $n$-convex for all $n\in\N$.

Assume now that $\{A_\gamma\mid\gamma\in\Gamma\}$ is a chain of $\F$-convex sets. 
Denote $A:=\bigcup_{\gamma\in\Gamma}A_\gamma$. Let $n\in\N$. To show the $n$-convexity of $A$,
let $x\in n^{-1}([n]A)$. Then $nx\in [n]A$, hence there exist $a_1,\dots, a_n\in A$ such that
$nx=a_1+\dots+a_n$. By the chain property, there exists a $\gamma\in\Gamma$ such that $a_1,\dots,a_n\in 
A_\gamma$. Hence $nx\in [n]A_\gamma$, which yields that 
\Eq{*}{
x\in n^{-1}([n]A_\gamma)\subseteq A_\gamma\subseteq A.
}
Thus $A$ is $n$-convex for all $n\in\F$, and hence it is $\F$-convex.

Now consider the case when $\{A_\gamma\mid\gamma\in\Gamma\}$ is a chain of $\F$-konvex sets. 
Denote $A:=\bigcup_{\gamma\in\Gamma}A_\gamma$. For $n\in\F$, to show the $n$-konvexity of $A$,
let $x\in [n]A$. By using the chain property again, it follows that there exists a $\gamma\in\Gamma$ 
such that $x\in [n]A_\gamma$. Hence, by the $n$-konvexity of $A_\gamma$, we have that $x\in 
nA_\gamma\subseteq nA$. Thus, $A$ is $n$-convex for all $n\in\F$, and hence it is $\F$-convex.

To show that the family of $\F$-konvex sets is not closed under the intersection in general, let $\F$ be the singleton $\{n\}$ with $n>1$ and consider $S:=\T=\R / \Z$. The complementary arcs $A=\big[0, \frac1n\big]$ and $B=\big[\frac1n, 0\big]$ are $n$-konvex. However, $n(A\cap B)=n\big\{0, \frac1n\big\}=\{0\}$ and $[n](A\cap B)=[n]\big\{0, \frac1n\big\}=\big\{0, \frac1n, \dots, \frac{n-1}n\big\}$, so the intersection $A\cap B$ is not $n$-konvex. 

Finally, let $A$ be an $\F$-convex and $B$ be an $\F$-konvex subset of $S$. Choose $n\in\F$. To prove 
that $A\cap B$ is $n$-konvex, let $x\in[n](A\cap B)$. Then, by the $n$-konvexity of $B$, we have that 
$x\in[n](A\cap B)\subseteq[n]B\subseteq nB$. Hence, there exists $b\in B$ such that $x=nb$. Therefore, 
the $n$-convexity of $A$ yields that $b\in n^{-1}(\{x\})\subseteq n^{-1}([n]A)\subseteq A$. Thus
$b\in A\cap B$ and consequently, $x=nb\in n(A\cap B)$, which proves the $n$-konvexity of $A\cap B$.
\end{proof} 

In view of the intersection property of $\F$-convex sets, for a nonempty subset $\F\subseteq\N$, we can 
define the $\F$-convex hull of an arbitrary subset $A$ of $S$ by
\Eq{*}{
  \conv_\F(A):=\bigcap\{C\mid A\subseteq C \subseteq S \mbox{ and } C \mbox{ is $\F$-convex}\}.
}
Clearly, $A\subseteq \conv_\F(A)$ and, by \prp{2}, $\conv_\F(A)$ is $\F$-convex. In particular, $A=\conv_\F(A)$ if and only if $A$ is $\F$-convex.  It is also obvious that, for $A\subseteq B\subseteq S$, we have $\conv_\F(A)\subseteq\conv_\F(B)$.

\Thm{CH}{Let $\F$ be a nonempty subset of $\N$. Then, for every $A\subseteq S$,
\Eq{CH}{
  A\subseteq \bigcup_{n\in\F}n^{-1}([n]A)\subseteq \conv_\F(A)
   \subseteq \bigcup_{n\in\langle\F\rangle}n^{-1}([n]A),
}
where $\langle\F\rangle$ denotes the multiplicative semigroup generated by $\F$.}

\begin{proof} The first inclusion in \eq{CH} is a consequence of \eq{I1}. The set $\conv_\F(A)$ is 
$\F$-convex and contains $A$, hence, for any $n\in\F$,
\Eq{*}{
   n^{-1}([n]A)\subseteq n^{-1}([n]\conv_\F(A)) \subseteq \conv_\F(A),
}
which proves the second inclusion in \eq{CH}. To prove the third one, it suffices to show 
that the set $C:=\bigcup_{n\in\langle\F\rangle}n^{-1}([n]A)$ is $\F$-convex. In fact, we will show that this 
set is actually $\langle\F\rangle$-convex. Let $k\in\langle\F\rangle$ and let
\Eq{*}{
 x\in k^{-1}\bigg([k]\bigg(\bigcup_{n\in\langle\F\rangle}n^{-1}([n]A)\bigg)\bigg).
}
Then, there exist $n_1,\dots,n_k\in \langle \F \rangle$ and $y_i\in n_i^{-1}([n_i]A)$ ($i\in\{1,\dots,k\}$) such that
\Eq{*}{
  kx=y_1+\cdots+y_k.
}
We have that $n_iy_i\in[n_i]A$, hence, with $n_0:=n_1\cdots n_k\in\langle\F\rangle$, it follows that
\Eq{*}{
  n_0y_i=\frac{n_0}{n_i}n_iy_i\in \frac{n_0}{n_i}([n_i]A)
   \subseteq \Big[\frac{n_0}{n_i}\Big]([n_i]A) = [n_0]A.
} 
Therefore,
\Eq{*}{
  n_0kx=n_0y_1+\cdots+n_0y_k\in [k]([n_0]A)=[kn_0]A.
}
This shows that
\Eq{*}{
  x\in (kn_0)^{-1}([kn_0]A))\subseteq \bigcup_{n\in\langle\F\rangle}n^{-1}([n]A)
}
since $kn_0\in\langle\F\rangle$. Thus, the $\langle\F\rangle$-convexity of the set $C$ is completed.
\end{proof}

The following result is an obvious consequence of \thm{CH}.

\Cor{CH}{Let $\F$ be a multiplicative subsemigroup of $\N$. Then, for every $A\subseteq S$,
\Eq{CH+}{
  \conv_\F(A)=\bigcup_{n\in\F}n^{-1}([n]A).
}}

Below we deal with the $\F$-convex hull of the union of finitely  many sets. The first observation follows from the definition of the $\F$-convex hull.

\Rem{r9}{Given a nonempty subset $\F$ of $\N$, for the $\F$-convex hull of the union of finitely many sets,
the following inclusion is always valid:
\Eq{DL}{
  \conv_\F(A_1\cup\cdots\cup A_k)
    \supseteq \bigcup\big\{\conv_\F(\{a_1,\dots,a_k\})\mid a_1\in \conv_\F(A_1),\,\dots,a_k\in \conv_\F(A_k)\big\}.
}}

The following result establishes an upper estimate for the $\F$-convex hull.

\Thm{DL}{Let $\F$ be a multiplicative subsemigroup of $\N$ and let $A_1,\dots,A_k\subseteq S$. 
If $A_1\subseteq B_1, \dots, A_k\subseteq B_k$ with some $\N$-konvex sets $B_1, \dots, B_k \subseteq S$, 
then
\Eq{12}{
  \conv_\F&(A_1\cup\cdots\cup A_k)\\
   &\subseteq \bigcup\big\{\conv_\F(\{b_1,\dots,b_k\})\mid b_1\in \conv_\N(A_1)\cap B_1,\,\dots,b_k\in \conv_\N(A_k)\cap B_k\big\}.
}}

\begin{proof} Clearly, we may assume that all the sets $A_1, \dots, A_k$ are nonempty. Assume that $A_1, \dots, A_k$  are covered by $\N$-konvex sets $B_1,\dots,B_k\subseteq S$, 
respectively. For the inclusion \eq{12}, in view of \cor{CH}, it is sufficient to show that, for all 
$n\in\F$,
\Eq{*}{
  n^{-1}\big([n]&(A_1\cup\cdots\cup A_k)\big)\\
   &\subseteq \bigcup\big\{\conv_\F(\{b_1,\dots,b_k\})\mid b_1\in \conv_\N(A_1)\cap B_1,\,\dots,b_k\in \conv_\N(A_k)\cap B_k\big\}.
}
To this aim, fix $n \in \F$ and take any $x \in n^{-1}\big([n]\big(A_1 \cup \dots \cup A_k \big)\big)$. Then 
\Eq{*}{
  nx=\sum^{n_1}_{j=1}x_{1,j}+ \dots + \sum^{n_k}_{j=1}x_{k,j},
} 
where $n_1, \dots, n_k$ are nonnegative integers summing up to $n$ and, for every $i \in \{1, \dots, k\}$,
\Eq{*}{
  x_{i,1}, \dots, x_{i,n_i}\in A_i\subseteq B_i
} 
(in fact we choose no $x_{i,j} \in A_i$ whenever $n_i=0$). For every $i \in \{1, \dots, k\}$, let $b_i$ be an 
arbitrary element of $B_i$ when $n_i=0$, and let $b_i$ be such a point of $B_i$ that 
\Eq{*}{
  x_{i,1}+ \dots + x_{i,n_i}=n_i b_i
}
otherwise; the last choice is possible due to the $n_i$-konvexity of $B_i$. 
Of course 
\Eq{*}{
  b_i\in n_i^{-1}(x_{i,1}+ \dots + x_{i,n_i})\subseteq n_i^{-1}([n_i]A_i)\subseteq \conv_\N(A_i).
}
Then we have
\Eq{*}{
  nx=n_1b_1+ \dots+ n_kb_k \in [n]\{b_1, \dots, b_k\},
}
hence
\Eq{*}{
  x\in n^{-1}\big([n]\{b_1, \dots, b_k \}\big) \subseteq \conv_\F(\{b_1,\dots,b_k\}),
}
which completes the proof of the inclusion \eq{12}.
\end{proof}

As an immediate consequence of  \thm{DL} and \rem{r9} we obtain the form of the $\N$-convex hull of the union of finitely  many sets covered by $\N$-konvex sets. 

\Cor{11}{If  $A_1, \dots, A_k \subseteq S$ are covered be $\N$-konvex subsets of $S$, then
\Eq{*}{
  \conv_\N(A_1\cup\cdots\cup A_k)
    = \bigcup\big\{\conv_\N(\{a_1,\dots,a_k\})\mid a_1\in \conv_\N(A_1),\,\dots,a_k\in \conv_\N(A_k)\big\}.
		}
}

As we mentioned in the introduction, singletons of $S$ may not be $\F$-convex (however, they are always 
$\F$-konvex). Therefore, it is of particular interest to describe the $\F$-convex hull of singletons.

\Prp{CH}{Let $\F$ be a multiplicative subsemigroup of $\N$. Then, for all $x,y\in S$, the $\F$-convex 
hulls of $x$ and $y$ are either disjoint, or equal.}

\begin{proof} By \cor{CH}, we have 
\Eq{CH0}{
 \conv_\F(\{x\})=\bigcup_{n\in\F}n^{-1}\{nx\} \qquad\mbox{and}\qquad 
 \conv_\F(\{y\})=\bigcup_{n\in\F}n^{-1}\{ny\}.
}
If $u$ is a common element of these sets, then there exist $n,m\in\F$ such that
\Eq{*}{
  nu=nx \qquad\mbox{and}\qquad mu=my.
}
Thus $nmx=nmy$. The set $\F$ being a multiplicative subsemigroup, $nm\in\F$, hence, in view of 
\eq{CH0}, we have
\Eq{*}{
  x\in\conv_\F(\{y\}) \qquad\mbox{and}\qquad y\in\conv_\F(\{x\}).
}
This yields that
\Eq{*}{
  \conv_\F(\{x\})\subseteq\conv_\F(\{y\}) \qquad\mbox{and}\qquad 
  \conv_\F(\{y\})\subseteq\conv_\F(\{x\}),
}
i.e., these $\F$-convex hull coincide.
\end{proof}

Based on the statement of this proposition, it follows that the family 
\Eq{*}{
  \big\{\conv_\F(\{x\})\mid x\in S\big\}
}
is a partition of $S$. The equivalence relation induced by this partition will be denoted by $x\sim_\F y$. 

\Prp{E}{Let $\F$ be a multiplicative subsemigroup of $\N$. Then, for any $x,y\in S$, the equivalence 
$x\sim_\F y$ holds if and only if there exists $n\in\F$ such that $nx=ny$. In addition, the semigroup 
operation of $S$ is compatible with $\sim_\F$, i.e., if $x_1\sim_\F x_2$ and $y_1\sim_\F y_2$, then 
$x_1+y_1\sim_\F x_2+y_2$.}

\begin{proof} For any $x,y\in S$ the equivalence 
$x\sim_\F y$ holds if and only if $x$ belongs to the $\F$-convex hull of $y$. This means, by \eq{CH0}, that
there exists $n\in\F$ such that $x\in n^{-1}\{ny\}$, which is the same as $nx=ny$. 

Now assume that $x_1\sim_\F x_2$ and $y_1\sim_\F y_2$. Then, there exist $n,m\in\F$ such that 
$nx_1=nx_2$ and $my_1=my_2$. Hence $nm(x_1+y_1)=nm(x_2+y_2)$, which proves that $x_1+y_1\sim_\F x_2+y_2$.
\end{proof}

In view of this proposition, the equivalence classes of the relation $\sim_\F$ form a commutative semigroup 
with respect to addition of subsets of $S$ which will be denoted by $\widetilde{S}^\F$. Analogously, 
for an element $x\in S$, the equivalence class containing $x$ will be denoted by $\widetilde{x}^\F$.

It is important to observe that in $\widetilde{S}^\F$ the $\F$-cancellation law holds.

\Prp{CL}{Let $\F$ be a multiplicative subsemigroup of $\N$. Then, for all $x,y\in S$ and $n\in\F$,
the equality $n\widetilde{x}^\F=n\widetilde{y}^\F$ implies $\widetilde{x}^\F=\widetilde{y}^\F$.}

\begin{proof} Indeed, if $n\widetilde{x}^\F=n\widetilde{y}^\F$, then $\widetilde{nx}^\F=\widetilde{ny}^\F$,
hence there exist $m\in\F$ such that $mnx=mny$. Since $mn$ belongs to $\F$, it follows that $x\sim_\F y$, and 
hence $\widetilde{x}^\F=\widetilde{y}^\F$.
\end{proof}

\section{A Stone-type theorem}

Given a subset $\F\subseteq\N$, two subsets $A,B$ of $S$ are called $\F$-{\it disjoint} if, 
for all $n \in \F$,
\Eq{*}{
[n]A \cap [n]B=\emptyset.
}
Obviously, the $\F$-disjointness of sets implies their disjointness, however, 
the converse may not be true. If $\F$ is a multiplicative subsemigroup of $\N$, 
then we have the following stronger statement.

\Prp{D}{Let $\F$ be a multiplicative subsemigroup of $\N$ and let $A,B\subseteq S$ be $\F$-disjoint subsets.
Then their $\F$-convex hulls are disjoint. Conversely, if one of the sets $A$ or $B$ is contained in an $\F$-konvex set and their $\F$-convex hulls are disjoint, then they are also $\F$-disjoint.}

\begin{proof} Assume that $x\in\conv_\F(A)\cap \conv_\F(B)$. Then,  by \cor{CH}, there exist $n,m\in\F$
such that $nx\in[n]A$ and $mx\in[m]B$. Hence, $nmx\in[nm]A\cap[nm]B$, which contradicts the 
$\F$-disjointness of $A$ and $B$.

To prove the reversed statement, assume that, $A$ is contained in an $\F$-konvex set $C\subseteq S$ and for some $n\in\F$, $x\in[n]A\cap[n]B$. Then, by the $n$-konvexity of $C$, $x\in [n]A\subseteq[n]C\subseteq nC$, i.e., for some $c\in C$, we have $x=nc$. Thus, $c\in n^{-1}([n]A)\cap n^{-1}([n]B)$, which contradicts the
disjointness of $\F$-convex hulls of $A$ and $B$.
\end{proof}

For complementary $\F$-disjoint sets, we have the following statement.

\Prp{l2}{Let $A$ and $B$ be $\F$-disjoint complementary subsets of $S$. Then $A$ and $B$ are $\F$-convex. 
If, in addition, $S$ is $\F$-konvex, then $A$ and $B$ are also $\F$-konvex.}

\begin{proof} Without loss of generality, we may assume that $A$ and $B$ are nonvoid. Let $n \in \F$ 
and $x \in S$ be such that $nx \in [n]A$. If $ x \not \in A$, then, by the complementarity of the sets $A$ 
and $B$, we get $x \in B$. Hence $nx \in nB \subseteq [n]B$, which contradicts the $\F$-disjointness of $A$ 
and $B$. Therefore, the condition $nx \in [n]A$ implies $x \in A$, and hence $A$ is 
$\F$-convex. The proof of the $\F$-convexity of $B$ is analogous.

The last assertion is a consequence of \prp{2}.
\end{proof}

The following lemma will play a crucial role in the proof of our Stone-type theorem.

\Lem{l1}{Let $A$ and $B$ be $\N$-disjoint subsets of $S$ and let $s \in S$. Then $A$ and $B\cup \{s\}$ or 
$A\cup \{s\}$ and $B$ are $\N$-disjoint.}

\begin{proof} Suppose on the contrary that 
\Eq{*}{
[n]A \cap [n](B\cup \{s\}) \not=\emptyset \quad {\rm and} \quad [m](A \cup \{s\})\cap [m]B \not= \emptyset
}
for some $n,m \in \N$. Then there exists $a_1, \ldots, a_n \in A$, $y_1, \ldots, y_n \in B\cup \{s\}$ and  $b_1, \ldots, b_m \in B$, $x_1, \ldots, x_m \in A\cup \{s\}$ such that 
\Eq{*}{
a_1+ \cdots+ a_n = y_1+ \cdots+ y_n
}
and
\Eq{eB}{
b_1+ \cdots+ b_m = x_1+ \cdots+ x_m.
}
Put $k:= \# \{i \in \{1, \ldots, m\}\mid x_i \in A\}$ and $l:=\# \{i\in \{1,\ldots, n\}\mid y_i \in 
B\}$. If $k=m$ then we would have $x_1, \ldots, x_m \in A$ and, consequently, $[m]B\cap [m]A \neq \emptyset$, 
which is impossible. Thus $0\leq k \leq m-1$ and, similarly, $0\leq l \leq n-1$. Without loss of generality, 
we may assume that $x_1, \ldots, x_k \in A$ and $y_1, \ldots, y_l \in B$.

Now consider the element
\Eq{*}{
u:=(n-l)(x_1+ \cdots+ x_k)+(m-k)(y_1+ \cdots+ y_l)+(m-k)(n-l)s.
}
Then, by \eq{eB}, 
\Eq{*}{
u&=(n-l)(x_1+ \cdots+ x_k +(m-k)s)+(m-k)(y_1+ \cdots+ y_l)\\
 &=(n-l)(x_1+ \cdots+ x_m)+(m-k)(y_1+ \cdots+ y_l)\\
 &=(n-l)(b_1+ \cdots+ b_m)+(m-k)(y_1+ \cdots+ y_l)\\
 &\in [(n-l)m+(m-k)l]B,
}
that is $u\in [nm-kl]B$. Similarly, we can prove that $u\in [nm-kl]A$, contrary to the 
$\N$-disjointness of $A$ and $B$.
\end{proof}

For the separation of $\N$-disjoint  sets we have the following Stone-type theorem.

\Thm{ST}{Let $A_0$ and $B_0$ be $\N$-disjoint subsets of $S$. 
Then there exist  $\N$-disjoint convex subsets  $A$ and $B$ such that
\Eq{ST}{
 A_0\subseteq A, \qquad B_0\subseteq B, \quad {\rm and} \quad A\cup B=S.
}
If, in addition, $S$ is konvex, then $A$ and $B$ are also konvex sets.}

\begin{proof}[1st Proof] By the Kuratowski--Zorn Lemma, there is a maximal element $(A,B)$ of the family 
\Eq{*}{
\{(A,B)\mid A_0\subseteq A\subseteq S,\, B_0\subseteq B\subseteq S,
 \,\,A\mbox{ and }B\mbox{ are $\N$-disjoint}\} 
}
partially ordered by the componentwise inclusion. Take any $s \in S$. Then, by \lem{l1}, $A$ and $B\cup\{s\}$ 
are $\N$-disjoint or $A\cup \{s\}$ and $B$  are $\N$-disjoint, and thus the maximality of $(A,B)$ yields 
that $s \in B$ or $s \in A$, respectively. This shows that $A\cup B=S$. Now the rest of the 
statement follows from \prp{l2}. 
\end{proof}

\begin{proof}[2nd Proof] We will deduce the statement using the Stone-type separation theorem of 
subsemigroups from the paper \cite{Pal89d}. 

Let $\bar{S}:=S\times\N$ and define the sets $\bar{A}_0$ and $\bar{B}_0$ by
\Eq{*}{
  \bar{A}_0:=\{(x,n)\mid n\in\N,\,x\in[n]A_0\} \qquad\mbox{and}\qquad
  \bar{B}_0:=\{(x,n)\mid n\in\N,\,x\in[n]B_0\}.
}
Then, obviously, $\bar{A}_0$ and $\bar{B}_0$ are subsemigroups of $\bar{S}$. The $\N$-disjointness
of $A_0$ and $B_0$ ensures that they are also disjoint. Thus, applying the Stone-type separation theorem of 
subsemigroups from \cite{Pal89d}, it follows that there exist disjoint subsemigroups $\bar{A}$ and $\bar{B}$ of $\bar{S}$ such that
\Eq{bar}{
  \bar{A}_0\subseteq \bar{A}, \qquad
  \bar{B}_0\subseteq \bar{B}, \qquad\mbox{and}\qquad
  \bar{A}\cup\bar{B}=\bar{S}.
}
Now define the sets $A,B\subseteq S$ by
\Eq{*}{
  A:=\{a\mid(a,1)\in\bar{A}\} \qquad\mbox{and}\qquad B:=\{b\mid(b,1)\in\bar{B}\}.
}
The inclusions in \eq{bar} show that \eq{ST} holds. The disjointness of $\bar{A}$ and $\bar{B}$ results
that $A$ and $B$ are $\N$-disjoint. Now the rest of the statement again follows from \prp{l2}. 
\end{proof}

%\nocite{JarLac09,JarLac10,Jar10b}
%\bibliography{publ,funcequ,t-conv}
%\bibliographystyle{plain}
%\end{document}

\end{document}